\documentclass[reqno,english,12pt]{amsart}
\pdfoutput=1                             

\usepackage{verbatim}
\usepackage[latin1]{inputenc}            
\usepackage[T1]{fontenc}                 
\usepackage[english]{babel}      
\usepackage{color}
\usepackage{fullpage}
\usepackage{tikz}
\usetikzlibrary{matrix}
\usetikzlibrary{arrows,chains,matrix,positioning,scopes}
\makeatletter
\tikzset{join/.code=\tikzset{after node path={%
\ifx\tikzchainprevious\pgfutil@empty\else(\tikzchainprevious)%
edge[every join]#1(\tikzchaincurrent)\fi}}}
\makeatother

\tikzset{>=stealth',every on chain/.append style={join},
         every join/.style={->}}

\usepackage[all]{xy}                     
\usepackage{mathtools}
\usepackage{euscript}

\newtheorem{definition}{Definition}

\newtheorem{theorem}{Theorem}
\newtheorem{lemma}[theorem]{Lemma}
\newtheorem{proposition}[theorem]{Proposition}
\newtheorem{corollary}[theorem]{Corollary} 

\newcommand{\R}{\mathbb{R}}
\newcommand{\N}{\mathbb{N}}

\newcommand{\Z}{\mathbb{Z}}
\newcommand{\Q}{\mathbb{Q}}

\def\co{\colon\thinspace}

\usepackage{amssymb}
\usepackage{yfonts}
\usepackage[numbers]{natbib}

\def\proof{\textbf{Proof.} }

\def\H{\mathbb{H} }

\def\tr{\text{tr}}

\newcommand{\Mod}[1]{\ (\mathrm{mod}\ #1)}

\begin{document}
\title{New examples from the jigsaw group construction}
\author{Carmen Galaz-Garc\'ia}
\begin{abstract}
A \textit{pseudomodular group} is a discrete subgroup $\Gamma \leq PGL(2,\Q)$ which is not commensurable with $PSL(2,\Z)$ and has cusp set precisely $\Q\cup\{\infty\}$. 
The existence of such groups was proved by Long and Reid. 
Later, Lou, Tan and Vo constructed two infinite families of non-commensurable pseudomodular groups which they called \textit{jigsaw groups}. 
In this paper we construct a new infinite family of non-commensurable pseudomodular groups obtained via this  jigsaw construction. 
We also find that infinitely many of the simplest jigsaw groups are not pseudomodular, providing a partial answer to questions posed by the aforementioned authors.  
\end{abstract}
\maketitle

\section{Introduction}

A \textit{Fuchsian group} $\Gamma$ is a discrete subgroup of $PSL(2, \R)$. 
Such a group acts properly discontinuously by fractional linear transformations on $\H^2$, the upper half-plane model of hyperbolic space. 
This $\Gamma$-action extends to the boundary at infinity $\partial_\infty \H^2 \equiv \R\cup\{\infty\}$.  
If an isometry $\gamma\in \Gamma$ has a fixed point on $\partial_\infty \H^2$ then it is either one of a pair of fixed points, in which case $\gamma$ is called \textit{hyperbolic}, or the fixed point is unique and we say $\gamma$ is \textit{parabolic}. 
The \textit{cusps} of a Fuchsian group $\Gamma$ are the points in $\partial_\infty\H^2$ fixed by parabolic elements of $\Gamma$. 
An example where the cusp set is easily calculated is when $\Gamma = PSL(2,\Z)$ is the modular group, in this case we obtain that  $\text{cusps}(PSL(2,\Z)) = \Q\cup\{\infty\}$. 

Two Fuchsian groups are \textit{commensurable} if they share a common subgroup that has finite index in both. 
It is known that commensurable Fuchsian groups  have the same cusp set. 
In \cite{LR} Long and Reid explore the converse question:  if $\Gamma_1$ and $\Gamma_2$ are finite covolume subgroups of $PSL(2, \R)$ with the same cusp set, are they commensurable? 
The answer was on the negative and in [2, theorem 1.2] they produced several examples of finite covolume Fuchsian groups with cusp set $\Q\cup\{\infty\}$ which are not commensurable with the modular group $PSL(2,\Z)$. 
This motivates the following definition:

\begin{definition}\em
A \textit{pseudomodular group} is a discrete subgroup $\Gamma \leq PGL(2,\Q)$ which is not commensurable with $PSL(2,\Z)$ and has cusp set precisely $\Q\cup\{\infty\}$.
\end{definition}

Subsequently Ayaka and Tan \cite{Ayaka} found another isolated example of a pseudomodular group and later Lou, Tan and Vo [3, theorem 1.2] constructed two infinite families of non-commensurable pseudomodular groups which they called \textit{jigsaw groups}. 
In this paper we examine a new infinite family of non-commensurable pseudomodular groups obtained via the jigsaw construction. 
We also find that infinitely many of the simplest jigsaw groups, called Weierstrass groups, are not pseudomodular. 

To describe the jigsaw construction from \cite{LTV} first let  $\Delta_n$, for $n\in \N$, be the ideal oriented triangle in $\H^2$ with vertices $\infty, -1$ and 0, and marked points 
\begin{equation}\label{markedPOINTS}
	x_1=-1+i, \ \ \ x_\frac{1}{n}=\frac{-n +i\sqrt{n}}{n+1}, \ \ \ x_n= i\sqrt{n}
\end{equation}
on the sides $ [\infty,-1]$, $[-1,0]$ and $[0,\infty]$ respectively. 
A \textit{tile of type} $n$ is any isometric transformation of $\Delta_n$, keeping track of the images of marked points. 
The sides of a tile of type $n$ will be called of \textit{type 1}, \textit{type} $\frac{1}{n}$ or \textit{type} $n$ according to which has the image of $x_1$, $x_\frac{1}{n}$ and $x_n$. 
Consider the $\pi$-rotations $\rho_i$ about the marked points $x_i$, represented here as elements of $PSL(2,\R)$: 
\begin{equation}\label{rotationGENERATORS}
	\rho_1=\begin{pmatrix}1 & 2 \\ -1 & -1 \end{pmatrix}, \ \ 
	\rho_{\frac{1}{n}}=\sqrt{n}\begin{pmatrix}1 & 1\\ \frac{-n-1}{n} & -1 \end{pmatrix}, \ \ 
	\rho_n= \frac{1}{\sqrt{n}}\begin{pmatrix} 0 & n \\ -1 & 0 \end{pmatrix}.
\end{equation}
\begin{definition}\label{Weir_groups}\em
The \textit{ $n$-th Weierstrass group} $W_n$ is the discrete group $W_n = \langle \rho_1, \rho_{\frac{1}{n}}, \rho_n \rangle.$
\end{definition}
For any $n\in \N$ the quotient surface $\H^2/W_n$ is an orbifold with a single cusp and three cone points of degree 2. 
Given the choice of marked points the element $\rho_1\rho_{\frac{1}{n}}\rho_n \in W_n$ is parabolic, so the tile $\Delta_n$ is \textit{balanced}. 
Since the vertices of $\Delta_n$ are in $\Q\cup\{\infty\}$, then $W_n\leq PSL(2,\Q)$ and all the vertices of the tiling of $\H$ generated by the action of $W_n$ on $\Delta_n$ are in $\Q\cup\{\infty\}$.  
In the notation of \cite{LTV} $\Delta_n = \Delta(1,1/n,n)$ and $W_n = \Gamma(1,1/n,n)$. 

By gluing different tiles together we can create groups that are more complex than the Weierstrass groups.
If we have two tiles $\Delta$ and $\Delta'$ with sides $s_1, s_2, s_3$ and $s'_1, s'_2, s'_3$, and marked points $x_1,x_2,x_3$ and $x'_1,x'_2,x'_3$ respectively, we say the sides $s_i$ and $s_j'$ match if both sides are of the same type. 
As explained in [3, definition 2.2], this means that if we glue $\Delta$ to $\Delta'$ along $s_i$ and $s_j'$ by identifying $x_i$ to $x_j'$, then the $\pi$-rotation about $x_i=x_j'$ will send $\Delta$ to $\Delta'$. 
In this way, by gluing finitely many tiles we obtain a triangulated ideal polygon with marked points on the interior and exterior sides of the triangulation, such a polygon is called a \textit{jigsaw}.

\begin{definition}\label{defjigsawgrp}\em
The \textit{jigsaw group} $\Gamma_J$ associated to a jigsaw $J$ is the Fuchsian group generated by the $\pi$-rotations about the marked points of the (exterior) sides of $J$. 
\end{definition}

As a convention we will require that the jigsaw $J$ used to define the jigsaw group $\Gamma_J$ has a tile $\Delta_n$ with vertices $\infty$, -1 and 0 in it.  
The balancing condition on each tile of the jigsaw $J$ ensures the quotient $\H^2/\Gamma_J$ is a complete orbifold with a single cusp and $N+2$ cone points of order 2, where $N$ is the number of tiles that make up the jigsaw. 
Then $\Gamma_J$ generates a tiling of $\H^2$ and $J$ is a fundamental domain of the action of the group. 

In [3, theorems 2.4 and 2.5] Lou, Tan and Vo examine jigsaw groups composed of tiles of types 1, 2 and 3. 
They prove that jigsaws composed only of tiles of types 1 and 2 have cusp set equal to $\Q\cup\{\infty\}$, 
and those that consist of a single tile of type 2 and $n$ tiles of type 1 are all pseudomodular and pairwise non-commensurable. 
On the other hand, they prove jigsaws made with tiles of type 1 and type 3 produce both an infinite family of pseudomodular groups and an infinite family of non-pseudomodular groups. 
Here we examine the groups generated by jigsaws made of tiles of types 1 and 4. 

\begin{figure}[h]
\begin{center}
\includegraphics[scale=1]{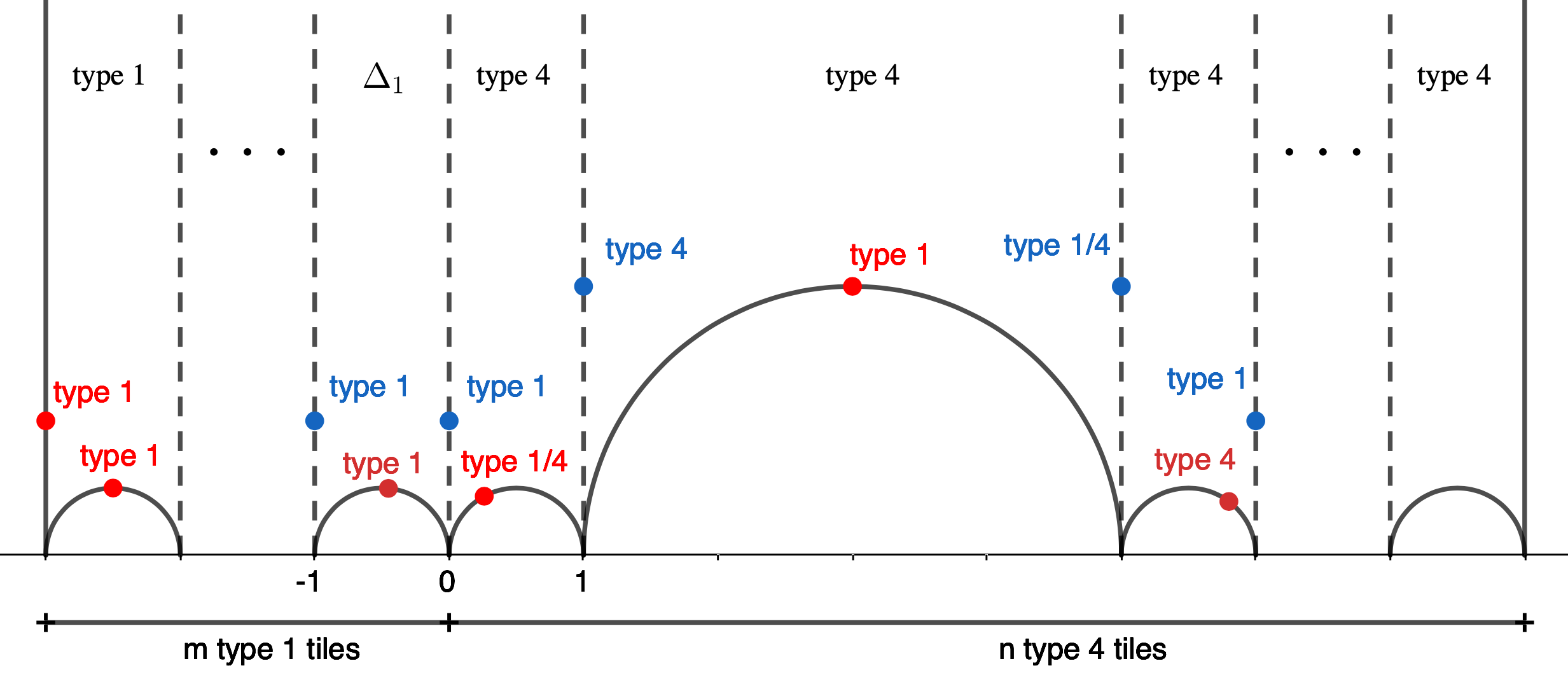} \ \ \ \ \ \ 
\end{center}
\caption{Jigsaw $J_{m,n}$}
\label{fig:Jmn}
\end{figure}

\begin{theorem}\label{pseudTHEOREM}
Let $J_{m,n}$ be the jigsaw formed by the $\Delta_1$ tile  followed by 
$m-1\geq0$ tiles of type 1 glued to the left and
$n\geq 1$ tiles of type 4 glued to the right of $\Delta_1$,
 so that all tiles in $J_{m,n}$ share $\infty$ as a common vertex (see figure \ref{fig:Jmn}). 
Then the associated jigsaw group $\Gamma_{m,n}$ has cusp set $\Q\cup\{\infty\}$. 
The infinite families $\Gamma_{1,n}$ and $\Gamma_{m,1}$ are pseudomodular and pairwise non-commensurable. 
\end{theorem}

To prove this we refine the process followed in \cite{LTV}. 
We first see the cusp set of these groups is $\Q\cup\{\infty\}$ by finding an explicit covering of $\R$ by killer intervals and then check they are non-commensurable by proving that each jigsaw group in the given families is non-arithmetic and equals its commensurator. 
By carefully analyzing the combinatorics of gluing tiles together we can extend the examples of pseudomodular groups to jigsaws with more than one tile of type $n>1$. 

In the final section of this paper we investigate whether there  are only finitely many pseudomodular Weierstrass groups $W_n$, a question posed by Lou, Tan and Vo in [3, section 9] which is a particular instance of the first open question posed by Long and Reid [2, section 6]. The following result provides a partial answer to these questions. 

\begin{theorem}\label{WnTHEOREM}
The groups $W_n$ with $n\geq 6$ and congruent to 0, 2 or 6 modulo 8 are not pseudomodular. 
\end{theorem}

For small values of $n\equiv 4 \Mod{8}$ we have found that $W_n$ is not pseudomodular. 
To construct these examples we have developed a computer program which tries to determine whether a given jigsaw group has cusp set equal to $\Q\cup\{\infty\}$ or contains a hyperbolic element fixing two rational points in $\partial_\infty\H^2 \equiv \R \cup \{\infty\}$. 
A survey of whether $W_n$ is pseudo-modular or not for $n\leq 28$ can be found at the end of section \ref{sectionWEIR}.

\section{Cusp set of the $\Gamma_{m,n}$ jigsaw groups.}

Let $\Gamma < PSL(2,\Q)$ be a Fuchsian group such that the quotient $\H^2/\Gamma$ has a single cusp.
Assume that $\infty$ is fixed by a parabolic element in $\Gamma$, so that the orbit of $\infty$ under the action of $\Gamma$ equals the cusp set of $\Gamma$. 
Since $\Gamma \leq PSL(2,\Q)$ then $\Gamma\cdot \infty \subseteq \Q\cup \{ \infty \}$. 
Therefore to prove $\text{cusps}(\Gamma) = \Q\cup\{\infty\}$ we only need to see that $\Q\subseteq \Gamma\cdot \infty$. 
To check this Long and Reid [2, example 1] introduced the following concept.

\begin{definition}\em
Let $p \in \Q$ be a cusp of $\Gamma$. 
A \textit{killer interval $I$ around $p$} is an interval $I \subset \R$ with $p\in I$
 for which there exists $\gamma \in \Gamma$ such that
 if $k\in I$ is a rational number, then the absolute value of the denominator of $\gamma(k)$ is strictly smaller than that of $k$. 
\end{definition}
If $\R$ can be covered by killer intervals  then for every $k\in \Q$ there will be a $\gamma \in \Gamma$ such that $\gamma(k) = \infty$.  
It is easy to see that every rational cusp of  $\Gamma$ has a killer interval around it. 
In detail, if $\gamma \in \Gamma < PSL(2,\Q)$ is parabolic then we can always find a matrix $g = \begin{pmatrix} a & b\\ c & d \end{pmatrix}\in PGL(2,\Q)$ such that $a,b,c,d\in \Z$, $\gcd(a, b, c, d)=1$ and both $\gamma$  and $g$ have the same action on $\H^2$.  
Then $\frac{a}{c}$ is a cusp of $\Gamma$ and $(\frac{a}{c} - \frac{1}{c}, \frac{a}{c} + \frac{1}{c})$ is a killer interval around it with associated map $\gamma$. 

\begin{definition}\em
Let $L = \begin{pmatrix}1 & 1 \\ 0 &1 \end{pmatrix}$ and suppose $\Gamma$ contains some power of $L$. Then
	 $$\ell(\Gamma) = \min \{ k\in \Z_{>0} \ |\   L^k \in \Gamma\}$$
is the \textit{fundamental length of $\Gamma$}.  A \textit{fundamental interval} for $\Gamma$ is any interval $[k,k+\ell(\Gamma)]$ with $k\in \Z$.
\end{definition}

If $I = [k, k+\ell(\Gamma)]$ is a fundamental interval for $\Gamma$ then every $x\in \R$ can be moved into $I$ by a power of $L^{\ell(\Gamma)}$. 
Translating by a a power of $L$ does not increase the denominator of a rational number. 
Then to prove $\Q\subseteq \Gamma\cdot \infty$ we just have to cover a fundamental interval of $\Gamma$ with killer intervals. 

Now let $J_{m,n}$ be the jigsaw formed by the $\Delta_1$ tile with vertices -1, 0 and $\infty$ followed by $m-1\geq0$ tiles of type 1 glued to its left and  $n\geq 1$ tiles of type 4 glued to its right, so that all tiles in $J_{m,n}$ share $\infty$ as a common vertex (see figure \ref{fig:Jmn}). 
Let $N=m+n$ and $v_0, v_1, \ldots, v_{N+1}$ be the cyclically ordered vertices of $J_{m,n}$, so that $v_0 = \infty, v_1 <v_2< \ldots <v_{N+1}$. 
For each $0\leq i \leq N$ let $x_i$ be the marked point on the side $[v_i,v_{i+1}]$ and $x_{N+1}$ be the marked point on $[v_{N+1}, v_0]$. 
Let $\Gamma_{m,n}$ be the jigsaw group associated to $J_{m,n}$.
If $\rho_i$ is the $\pi$-rotation around $x_i$ then $\Gamma_{m,n} = \langle \rho_0, \rho_1, \ldots, \rho_{N+1}\rangle$. 
Clearly the vertices $v_1, \ldots, v_{N+1}$ are in the orbit of $v_0$ so $\H^2/\Gamma_{m,n}$ has a single cusp. 
Proposition 4.5 in [3] proves that 
\begin{equation}\label{ellGamma}
	\rho_{N+1}\rho_N\ldots \rho_0 = \begin{pmatrix}1 & \ell(\Gamma_{m,n}) \\ 0 & 1 \end{pmatrix}
\end{equation}
and $\ell(\Gamma_{m,n}) = 3m + 6n$.
Then $v_0=\infty$ is fixed by a parabolic element of $\Gamma_{m,n}$. 
This implies that every vertex of a tile in the triangulation of $\H^2$ induced by $J_{m,n}$ is a cusp of $\Gamma_{m,n}$.


In the following let $\Delta(a,b,c)$ be the ideal triangle with vertices $a,b,c$ and sides $[a,b]$, $[b,c]$ and $[c,a]$.
Denote the $\pi$-rotation about a point $(x,y)\in \H^2$ by $R_{x,y}$. 
When $x, y \in \Q$ it is possible to represent $R_{x,y}$ as a matrix $\begin{pmatrix}a& b \\ c& d\end{pmatrix}\in PGL(2,\Q)$ with $\gcd(a,b,c,d)=1$, this will allow us to calculate lengths of killer intervals.

In the triangulation of $\H^2$ produced by a jigsaw a \textit{vertical tile}  is one that has $\infty$ as a vertex. 
A \textit{vertical side} of the triangulation induced by a jigsaw on $\H^2$ is one that has $\infty$ as an endpoint, it can be interior or exterior. 

\begin{proposition}[\text{4.3 in [3]}] \label{vertical_tiles}
Let $T=\Delta(\infty, x_1, x_2)$ be a vertical tile of type 4 in the triangulation of $\H^2$ produced by a jigsaw. 
Let the sides of $T$ be $e_1=[\infty, x_1]$, $e_2=[x_1, x_2]$ and $e_3=[x_2,\infty]$ , so that $e_i$ has type $k_i$ and marked point $p_i$. 
Then there are three possible configurations for $T$:
\begin{itemize}
\item if $k_1 = 1$ then $k_2= \frac{1}{4}$, $k_3 = 4$ and $x_2 = x_1 +1$. The marked points are $p_1 = (x_1,1)$, $p_2 = (x_1 + \frac{1}{5}, \frac{2}{5})$ and $p_3=(x_1+1,2)$. The vertical tile to the right of $T$ has type 4. 
\item if $k_1 = 4$ then $k_2= 1$, $k_3 = \frac{1}{4}$ and $x_2 = x_1 +4$. The marked points are $p_1 = (x_1,2)$, $p_2 = (x_1+2,2)$ and $p_3=(x_1+4,2)$. The vertical tiles to the right and left of $T$ have type 4. 
\item if $k_1 = \frac{1}{4}$ then $k_2= 4$, $k_3 = 1$ and $x_2 = x_1 +1$. The marked points are $p_1 = (x_1,2)$, $p_2 = (x_1+\frac{4}{5},\frac{2}{5})$ and $p_3=(x_1+1,1)$. The vertical tile to the left of $T$ has type 4. 
\end{itemize}
All vertical tiles of type 1 are of the form $\Delta(\infty, x, x+1)$ with $m\in \Z$. The marked points on the sides $[\infty, x]$, $[x, x+1]$ and $[x+1,\infty]$ are $(x,1)$, $(x+\frac{1}{2}, \frac{1}{2})$ and $(x+1,1)$ respectively. 
\end{proposition}
In all figures a solid line indicates an exterior side of a tile and a dashed line indicates an exterior side. 
Dotted lines indicate sides that could be either interior or exterior. 


\begin{proposition}\label{cuspsareQ}
Let $J_{m,n}$ be a jigsaw as in theorem \ref{pseudTHEOREM} and $\Gamma$ its associated jigsaw group. 
Then $\text{cusps}(\Gamma)=\Q \cup \{\infty\}$.
\end{proposition}
\proof 
We will prove there is a covering of $\R$ by killer intervals of cusps of $\Gamma$. 
Consider the triangulation of $\H^2$ generated by the action of $\Gamma$ on the triangulated jigsaw $J_{m,n}$. 
Since all tiles in $J_{m,n}$ are of type 1 or 4, proposition 4.3 in \cite{LTV} implies that the vertices of a vertical tile that lie on $\R$ are integers at distance 1 or 4 from each other. 
Then $\R$ can be divided into consecutive intervals of lengths one and four, with each endpoint being an integer. 

If $v$ is an endpoint of a vertical side with $v\neq \infty$, 
then $v$ is a cusp of the jigsaw group and by proposition 4.6 in \cite{LTV} the killer interval around $v$ is $(v-1,v+1)$. 
Then to cover $\R$ with killer intervals  it will be enough to cover the gaps of length 4 between cusps.

By proposition \ref{vertical_tiles} a vertical tile $T_0$ with vertices  $m$ and $m+4$ has to be a tile of type 4 where the side $[m,m+4]$ is of type 1. 
Without loss of generality we may translate this tile and assume $T_0 = \Delta(\infty, 0, 4)$, its marked points are $(0,2)$, $(2,2)$ and $(4,2)$. 
Let $T_1 = R_{2,2}(T_0)$, so $T_1$ is adjacent to $T_0$ along the side $[0,4]$ and has vertices 4, 0 and 2. 
\\

\begin{figure}[h]
\begin{center}
\includegraphics[scale=2]{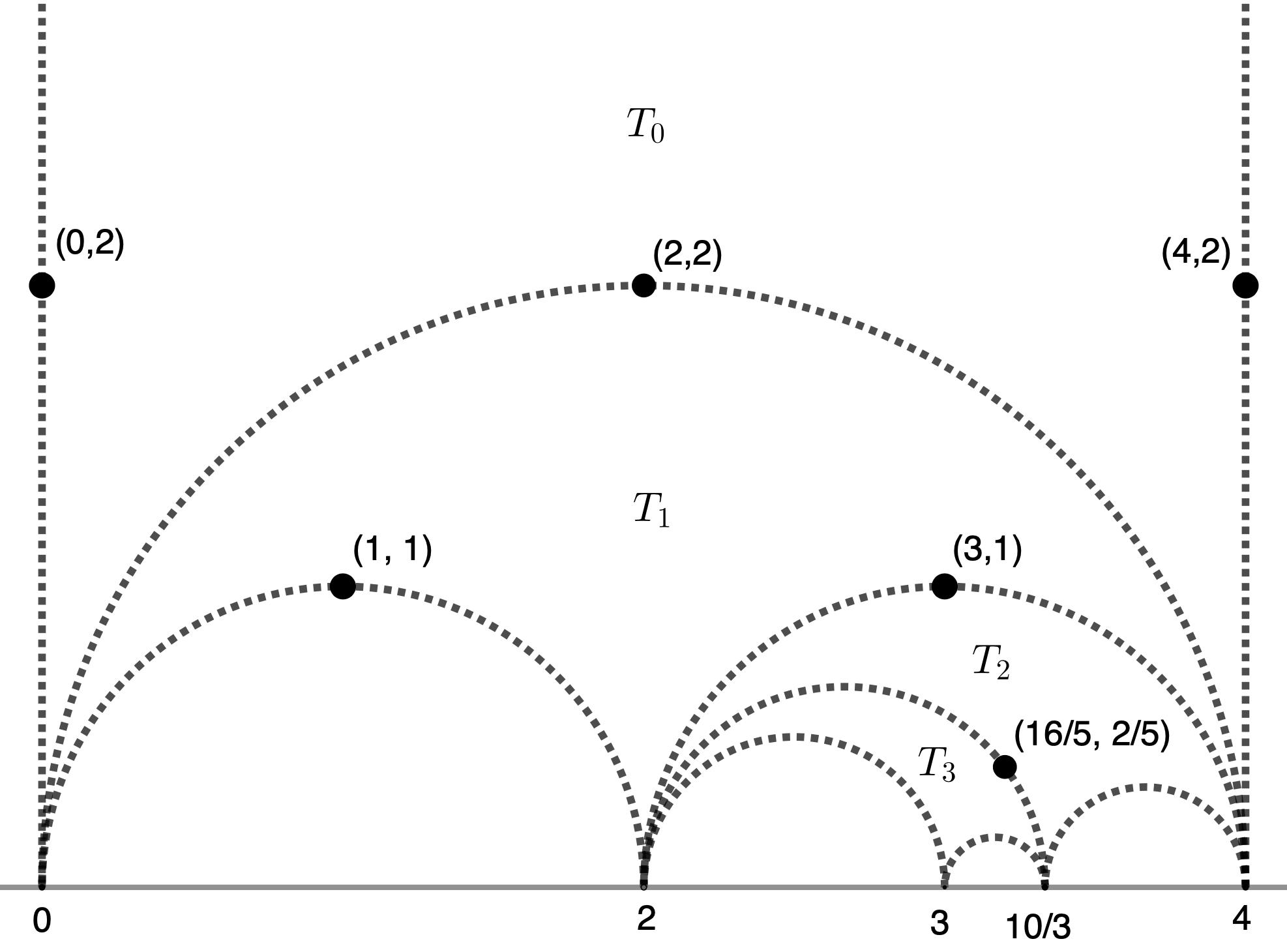} \ \ \ \ \ \ 
\end{center}
\caption{Tile $T_1$ is of type 4}
\label{fig:case2twotype1}
\end{figure}

\textbf{Case 1: $T_1$ is a tile of type 4.} 
For this case see figure \ref{fig:case2twotype1}. 
Since $0$ and $4$ are endpoints of a vertical side the killer intervals around these cusps are $(-1,1)$ and $(3,5)$. 
The tile $T_1$ is type 4, so by proposition 4.7 in \cite{LTV} the killer interval around $2$  is $(1,3)$. 
Then to cover the interval $[0,4]$ it will be enough to check that $1$ and $3$ are cusps of $\Gamma$. 
We will use the following matrices for calculations:
	$$R_{2,2} = \begin{pmatrix}2 & -8\\1  & -2\end{pmatrix}, \  R_{3,1}=\begin{pmatrix} -3 & 10 \\ -1& 3\end{pmatrix}, \ R_{\frac{16}{5},\frac{2}{5}} =\begin{pmatrix}  -16& 52 \\ -5 & 16\end{pmatrix}.$$

Since $T_1 = R_{2,2}(T_0)$, then the side $[2,4]$ of $T_1$ is type 4 with marked point $R_{2,2}(0,2)=(3,1)$. 
The tile adjacent to $T_1$ along $[2,4]$ is $T_2 = R_{3,1}(T_1) =  \Delta(4,2,\frac{10}{3})$, it is of type 4 as well. 
The side $[2,\frac{10}{3}]$ has type 1 with marked point $R_{3,1}(2,2)=(\frac{16}{5},\frac{2}{5})$. 
Finally, consider the tile $T_3$ that is adjacent to $T_2$ along $[2,\frac{10}{3}]$.
We have that $T_3 = R_{\frac{16}{5},\frac{2}{5}} (T_2) = \Delta(2,\frac{10}{3},3)$. 
This proves 3 is a vertex of a tile in the triangulation and therefore a cusp of $\Gamma$. 
By a similar argument we can prove 1 is a cusp of $\Gamma$. 
Notice this case also covers all jigsaws of the form $J_{0,n}$.

\textbf{Case 2: $T_1$ is a tile of type 1.} 
Since $0$ and $4$ are endpoints of vertical sides, the killer intervals around them still are $(-1,1)$ and $(3,5)$. 
The tile $T_1$ now has type 1, so by proposition 4.7 \cite{LTV} the killer interval around 2 is $(\frac{3}{2}, \frac{5}{2})$. 
We will see that $(1,\frac{5}{3})$ and $(\frac{7}{3},3)$ are killer intervals for $\frac{4}{3}$ and $\frac{8}{3}$ respectively. 
Then the killer intervals for  0, $\frac{4}{3}$, 2, $\frac{8}{3}$ and 4 will cover $[0,1]\setminus\{1,3\}$. 
To finish it will only be necessary to check that 1 and 3 are cusps of $\Gamma$.

\begin{figure}[h]
\begin{center}
\includegraphics[scale=1.4]{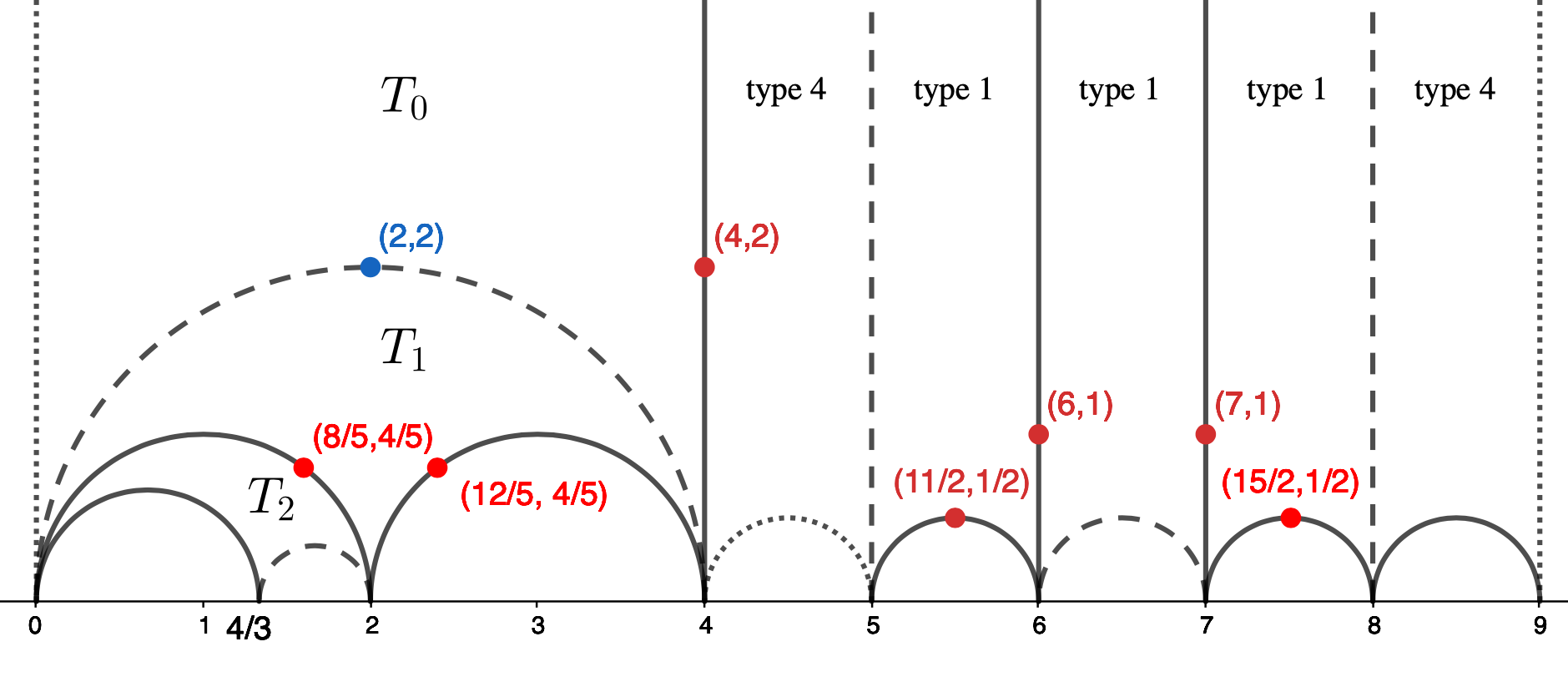} \ \ \ \ \ \ 
\end{center}
\caption{Vertical tiles for a jigsaw $J_{1,n}$}
\label{fig:jigsaw1tilestype1}
\end{figure}

\begin{figure}[h]
\begin{center}
\includegraphics[scale=1.5]{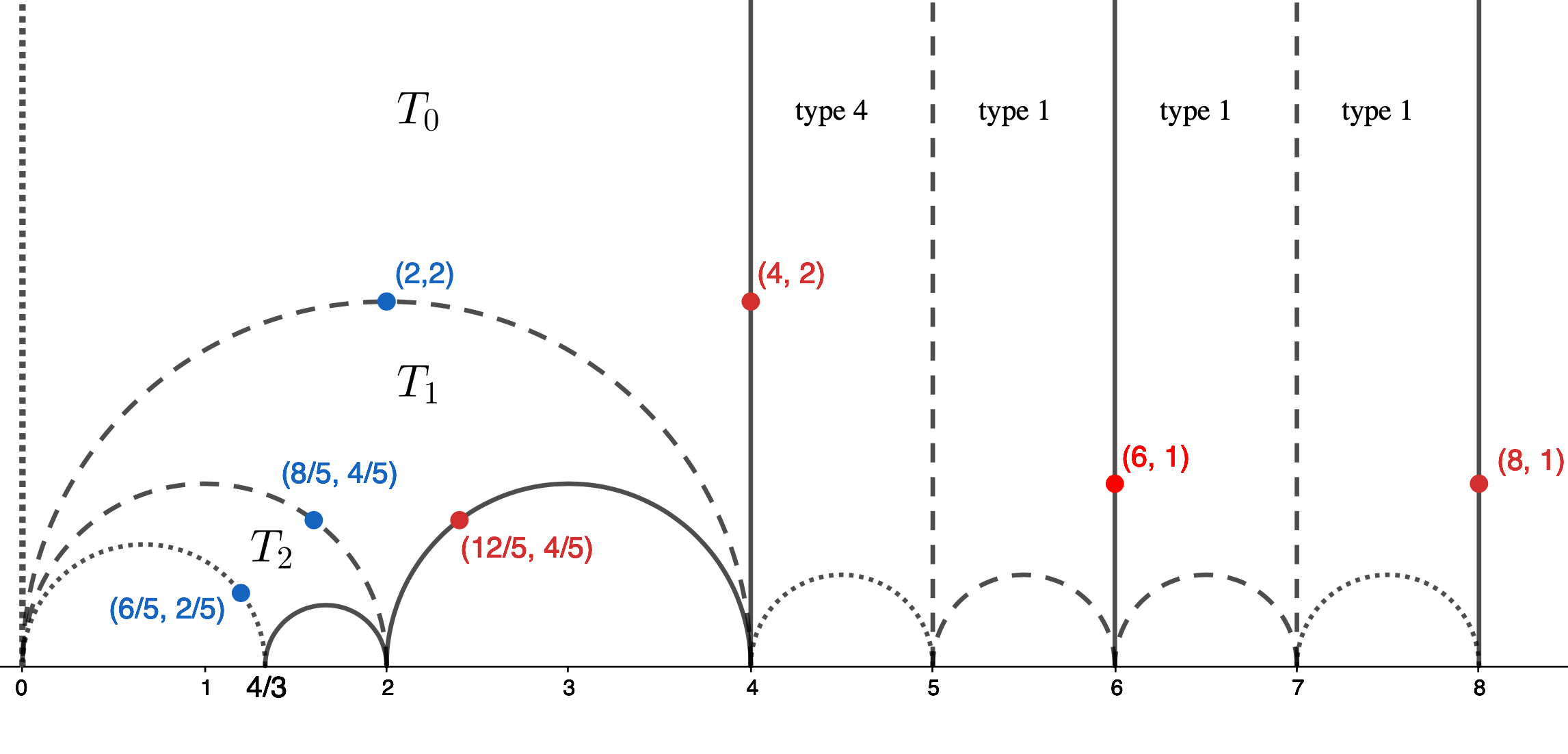} \ \ \ \ \ \ 
\end{center}
\caption{Vertical tiles for a jigsaw $J_{m,n}$ with $m\geq 2$}
\label{fig:jigsaw2tilestype1}
\end{figure}

If the side $[0,4]$ of $T_0$ was exterior, then by rotating around the marked point $(2,2)$ we would get that $T_1$ is also of type 4. 
Thus it must be that $[0,4]$ is an interior side. 
Then $T_0$ is in the $\Gamma$-orbit of the unique tile $T_0'$ of type 4 in the initial jigsaw $J_{m,n}$ that shares an interior side with a tile of type 1. 
Since $T_0'$ has an exterior side of type $\frac{1}{4}$ (see figure \ref{fig:Jmn}) then in $T_0$ the side $[4,\infty]$, which has type $\frac{1}{4}$, must be exterior too. 
The side $[0,\infty]$ of $T_0$ is only exterior when $n=1$. 
For a jigsaw $J_{1,n}$  the tiling follows the pattern shown in figure \ref{fig:jigsaw1tilestype1} and for a jigsaw $J_{m,n}$ with $m\geq 2$ the tiling is as in figure \ref{fig:jigsaw2tilestype1}.

\vspace{ 0.5 cm}

$\bullet$ 3 is a cusp. 
Notice that $[4,\infty]$ is an exterior side with marked point $(4,2)$ and 8 is a vertex of the tiling for every $J_{m,n}$. 
Then $R_{4,2}(3) = 8$ implies 3 is a cusp of $\Gamma$. 
To find an element of $\Gamma$ that sends $\infty$ to $3$ we will have to consider two cases. 
For a jigsaw $J_{1,n}$ the sides $[\infty,7]$ and $[7,8]$ are exterior of type 1 with marked points $(7,1)$ and $(\frac{15}{2},\frac{1}{2})$ respectively. 
Rotating around these points we get that $S = R_{4,2}R_{\frac{15}{2},\frac{1}{2}}R_{7,1} \in \Gamma$ sends $\infty$ to 3. 
For  $J_{m,n}$ with $m\geq 2$  the side $[8,\infty]$ is exterior with marked point $(8,1)$. 
Therefore $S' = R_{4,2}R_{8,1} \in \Gamma$ sends $\infty$ to 3. 
It can be calculated that
	$$S  = \begin{pmatrix} 12 & -64 \\ 4 & -21 \end{pmatrix}\  \text{ and }\  S'= \begin{pmatrix} 12 & -100 \\ 4 & -33 \end{pmatrix},$$
so in both cases we obtain that $(3 - \frac{1}{4}, 3 + \frac{1}{4})$ is a killer interval around 3.

$\bullet$ $\frac{4}{3}$ is a cusp. 
The side $[2,4]$ is exterior with marked point $(\frac{12}{5},\frac{4}{5})$ in every $J_{m,n}$, so $R_{\frac{12}{5},\frac{4}{5}} \in \Gamma$. 
Since 3 is a cusp of $\Gamma$, then $\frac{4}{3} = R_{\frac{12}{5},\frac{4}{5}}(3)$ is too.
For $J_{1,n}$ we get that $R_{\frac{12}{5},\frac{4}{5}}S \in \Gamma$ sends $\infty$ to $\frac{4}{3}$. For $J_{m,n}$ with $m\geq2$ we see that $R_{\frac{12}{5},\frac{4}{5}} S'\in \Gamma$ does the same. By making
	$R_{\frac{12}{5},\frac{4}{5}}=\begin{pmatrix} 12 & -32 \\ 5 & -12  \end{pmatrix}$
we have that
	$$R_{\frac{12}{5},\frac{4}{5}}S= \begin{pmatrix} 4 & -24 \\ 3 & -17\end{pmatrix}  \ \text{ and } \  
	R_{\frac{12}{5},\frac{4}{5}} S'  =\begin{pmatrix} 4 & -36 \\ 3 & -26\end{pmatrix}.$$
This shows the killer interval around $\frac{4}{3}$ is $(1,\frac{5}{3})$. 

 
$\bullet$ 1 is a cusp. 
The marked point on the side $[0,2]$ of $T_1$ is $(\frac{8}{5},\frac{4}{5})$, then $T_2=R_{\frac{8}{5},\frac{4}{5}}(T_1) = \Delta(0,\frac{4}{3},2)$  is in the triangulation. 
For every $J_{m,n}$ the tile $T_2$ is of type 1, 
so the marked point on the side $[0,\frac{4}{3}]$ is $(\frac{6}{5}, \frac{2}{5})= R_{\frac{8}{5},\frac{4}{5}}((2,2))$. 
Therefore $1 = R_{\frac{6}{5}, \frac{2}{5}}(2)$ is a cusp.


$\bullet$ $\frac{8}{3}$ is a cusp. 
Since $R_{\frac{12}{5},\frac{4}{5}}\in \Gamma$ for all $J_{m,n}$, then $\frac{8}{3} = R_{\frac{12}{5},\frac{4}{5}}(0)$ is a vertex of the tiling and  a cusp of $\Gamma$.
To find an element in $\Gamma$ that sends $\infty$ to $\frac{8}{3}$ recall that $T_0$ is in the $\Gamma$-orbit of the unique tile of type 4 in $J_{m,n}$ that has an interior side adjacent to a tile of type one. 
Then there must be an $n\in \N$ and $G\in\Gamma$ so that the tile $T=\Delta(\infty, n, n+1)$ is in the triangulation, has sides $[\infty, n]$, $[n,n+1]$ and $[n+1,\infty]$ of types 1, $\frac{1}{4}$ and 4 respectively, and $G(T)=T_0$. In particular we have that $G(\infty) = 0$. 
A direct calculation shows we can write $G = \begin{pmatrix}0& 4 \\ -1 & n+1 \end{pmatrix}$. Thus
	$R_{\frac{12}{5},\frac{4}{5}} G = \begin{pmatrix}8 & 4-8n \\ 3 & 2-3n \end{pmatrix} \in \Gamma$  sends $\infty$ to $\frac{8}{3}$. 
This shows the killer interval around $\frac{8}{3}$ is $(\frac{7}{3}, 3)$.
\vspace{0.5 cm}


\section{Non-commensurability of the $\Gamma_{1,n}$ and $\Gamma_{m,1}$ jigsaw groups.}

The \textit{commensurator} of a subgroup $\Gamma$ of $PSL(2,\R)$ is the subgroup 
	$$\text{Comm}(\Gamma)= \{g\in PSL(2,\R) \ | \ g\Gamma g ^{-1} \text{ commensurable with } \Gamma \}.$$
It is a theorem by Margulis \cite{Margulis} that if $\Gamma$ is non-arithmetic then $\text{Comm}(\Gamma)$ is the unique maximal element (with respect to subgroup inclusion) in the commensurability class of $\Gamma$. 
Following sections 7 and 8 in [3], to see that jigsaw groups $\Gamma$ of the form $\Gamma_{1,n}$ and $\Gamma_{m,1}$ are pairwise non-commensurable we will check that each $\Gamma$ is non-arithmetic and $\Gamma = \text{Comm}(\Gamma)$. 
To prove the latter we analyze the location of tangency points on the maximal horocycle of the orbifold $\H^2/\Gamma$.


\subsection{Non-arithmeticity.}
By Takeuchi \cite{Takeuchi} if a non-compact Fuchsian group $\Gamma \leq PSL(2,\R)$ of finite covolume, with no elements of order 2 and with invariant trace field $\Q$ is arithmetic, then $\tr(\gamma^2) \in \Z$ for all $\gamma \in \Gamma$. 
Since $\tr(\gamma^2) = (\tr\gamma)^2 -2$ it is enough to see whether $(\tr\gamma)^2\in \Z$. 

Let $J$ be a jigsaw as in theorem \ref{pseudTHEOREM} with associated jigsaw group $\Gamma$, and let $\rho_0, \ldots, \rho_{N+1}$ be the generators of $\Gamma$ as in (\ref{ellGamma}). 
We will see  the subgroup of index two $\Gamma^{(2)}$ consisting of all elements of $\Gamma$ with even word length is non-arithmetic, and therefore $\Gamma$ is non-arithmetic too. 
Notice the group $\Gamma^{(2)}$ still has finite covolume, a fundamental domain for $\Gamma^{(2)}$ is $J\cup\rho_0(J)$. 


\begin{proposition} \label{nonarith}
Let $J_{m,n}$ be a jigsaw as in theorem \ref{pseudTHEOREM} and $\Gamma$ its associated jigsaw group. 
Then $\Gamma$ is non-arithmetic.
\end{proposition}

\proof  It is enough to see that there exists $\gamma\in\Gamma^{(2)}$ such that $\tr(\gamma)^2 \notin \Z$.
Let $e_j =[x_j, \infty]$ and $e_k = [x_k, \infty]$ be exterior vertical sides  in the tiling of $\H^2$ induced by $J_{m,n}$. 
Assume that $e_j$ is of type 4 and $e_k$ is of type 1, so their  marked points are $(x_j, 2)$ and $(x_k,1)$ respectively. 
Since $e_j$ and $e_k$ are exterior sides the $\pi$-rotations
	$$R_{x_j,2}= \frac{1}{2}\begin{pmatrix} x_j & -(x_j^2+4) \\ 1 & -x_j \end{pmatrix} \text{ and } =R_{x_k,1}= \begin{pmatrix} x_k & -(x_k^2+1) \\ 1 & -x_k \end{pmatrix}$$
are elements of $\Gamma$.
We have that $(\tr (R_{x_j,2} R_{x_k,2}))^2 = \frac{1}{4}(-(x_k-x_j)^2 -5)^2$, so
\begin{eqnarray*}
	(\tr (R_{x_j,2} R_{x_k,2}))^2  \in \Z &\Leftrightarrow& (-(x_k-x_j)^2 -5)^2 \equiv 0\Mod{4} \\
					&\Leftrightarrow& -(x_k-x_j)^2 -5 \equiv 0\Mod{2}\\
					&\Leftrightarrow& (x_k-x_j)^2 \equiv 1\  \Mod{2}\\
					&\Leftrightarrow& x_k-x_j \equiv 1\ \Mod{2}.
\end{eqnarray*} 
Then if $\Gamma^{(2)}$ is arithmetic the distance between the real vertex of a vertical side of type 1 and the real vertex of a vertical side of type 4 must be  odd.
However, the jigsaw $J_{m,n}$ has a tile  $T$ of type 1 with two exterior sides, so there is a tile in the orbit $\Gamma \cdot T$ where both exterior type 1 sides are vertical and at distance one from each other. 
Therefore one of these consecutive exterior vertical sides of type 1 will be at even distance from a vertical side of type 4. 

\vspace{0.5 cm}


\subsection{Tangency points of maximal horocycle.}
Let $\Gamma$ be the jigsaw group associated to a jigsaw $J = J_{m,n}$ as in theorem \ref{pseudTHEOREM}. 
Then the orbifold $\mathcal{O} = \H^2/\Gamma$ has $N=m+n$ cone points of order 2, a cusp and finite volume.
Let $\pi \co \H^2 \to \mathcal{O}$ be the corresponding quotient map. 
The lift of the cone points of $\mathcal{O}$ to $\H^2$ is the set of all marked points on exterior sides in the tiling $\Gamma\cdot J$ of $\H^2$. 
Since $J_{m,n}$ only has tiles of type 1 and type 4, by proposition \ref{vertical_tiles} all the marked points in the tiling are on or below the line $y=2$. 

Recall that a \textit{horocycle in} $\H^2$ centered at $\xi \in \partial_\infty \H^2 \equiv \R\cup\{\infty\}$ is a curve $\alpha\setminus \{\xi\} \subset \H^2$ where $\alpha$ is a Euclidean circle tangent to $\R$ at $\xi$, if $\xi\in \R$, or $\alpha$ is a line parallel to the $x$-axis if $\xi =\infty$. 
A curve $C$ in $\mathcal{O}$ is a horocycle if $C$ is the image under $\pi$ of a horocycle in $\H^2$ and does not self-cross.
For $t>0$ let $\alpha_t$ be the line $y= t$. 
When $t>2$ the horocycle $\pi(\alpha_t)$ loops once around the cusp of $\mathcal{O}$ without self-intersecting and the length of $\pi(\alpha_t)$ goes to 0 as $t$ goes to $\infty$. 
The maximal horocycle in $\mathcal{O}$ is then $C = \pi(\alpha_2)$.
The curve $C$ is tangent to itself at the cone points of $\mathcal{O}$ that are projections of marked points of exterior sides in $\H^2$ with  $y$-coordinate equal to 2. 
The lift $\tilde{C}=\pi^{-1}(C)$ to $\H^2$ is formed by the horizontal horocycle $\alpha_2$, horocycles of radius 1 which are tangent to $\alpha_2$, and smaller horocycles based at the other cusps which are disjoint from $\alpha_2$. 

To prove that $ \Gamma_{m,n} = \text{Comm}(\Gamma_{m,n}) $ when $m=1$ or $n=1$ we will analyze the location of tangency points of $\tilde{C}$ along  $\alpha_2$. 
This will be used to see that the orbifold $\H^2/\text{Comm}(\Gamma_{m,n})$ cannot be finitely covered by $\H^2/\Gamma_{m,n}$. 
Recall that the horizontal translation in $\H^2$ by $\ell(\Gamma_{m,n}) = 3m + 6n$ is the smallest horizontal translation that is an element of $\Gamma$. 

\begin{lemma} \label{tanpoints1}
Let $J=J_{m,1}$ with $m\geq 1$ and associated jigsaw group $\Gamma$. 
Let $T$ be a horizontal translation by less than $\ell(\Gamma)$. 
Then there is a pair of tangency points  $p_1,p_2$ of  $\tilde{C}$ such that, 
if $L\neq Id$ is a horizontal translation by less than $\ell(\Gamma)$, then $L(p_1)$ and $L(p_2)$ are no longer tangency points of $\tilde{C}$.

\end{lemma}
\proof 
The tile $ E = \Delta(\infty,0,1)$ is the unique tile of type 4 in $J$. 
Its sides $[0,1]$ and $[1,\infty]$ are exterior of type $\frac{1}{4}$ and 4 respectively, the  marked point on $[1,\infty]$ is $p_1=(1,2)$. 
Then $E' = R_{1,2}(E) = \Delta(\infty, 1, 5)$ is of type 4 and has $[5,\infty]$ as an exterior side with marked point $p_2=(5,2)$.
The next tile $E'' = R_{5,2}(E') = \Delta(\infty,5,6)$ is of type 4 but now $[6,\infty]$ is an interior side and therefore there is no (exterior) marked point on it. 
Since $E$ is the unique tile of type 4 in $J$, the tiles $E$, $E'$ and $E''$ are the only vertical tiles of type 4 with vertices on the fundamental interval $I=[0,3m+6]$. 
Therefore $p_1$ and $p_2$ are the only two tangency points of $\tilde{C}$ at height 2 on $I\times[0,\infty)$. 
If $L$ is a translation by $0<k<\ell(\Gamma)$ where $L(p_1)$ is a tangency point of $\tilde{C}$, then it must be that $k=4$ and $L(p_1)=p_2$. 
But $L(p_2)$ is on a vertical type 1 side, so it cannot be a tangency point.

\vspace{0.5 cm}

To prove a similar result for jigsaws $J_{1,n}$ with $n>1$ we will need not a pair but a triple of tangency points on $\alpha_2$. 
To find these we examine patterns of consecutive vertical tiles.

\begin{definition}\em 
The \textit{width} of a vertical tile is the distance between its vertices on the $x$-axis. 
For $i=1, \ldots, k$ let $T_i = \Delta(\infty, x_i, x_{i+1})$  be a vertical tile with $x_i < x_{i+1}$ and width $w_i$. 
The \textit{width pattern} of the consecutive tiles $T_1, \ldots, T_m$ is the tuple $(w_1, \ldots, w_m)$. 
\end{definition}

By proposition \ref{vertical_tiles}, tiles of type 1 always have width 1 and tiles of type 4 have width either 1 or 4. 
In the proof of \ref{tanpoints2} we will also use \textit{half tiles}, these are translations of either $\Delta(\infty,0,4)\cap([0,2]\times \R)$ or  $\Delta(\infty,0,4)\cap([2,4]\times \R)$. 
Half tiles have width 2 and can be included in a list of adjacent vertical tiles to generate a width pattern. 
The width pattern $(2,2)$ is allowed to indicate two halves of the same tile of type 4 and width 4. 
We will need half tiles to account for marked points on the non-vertical sides of tiles of type 4. 
Finally, notice that not any tuple with coordinates 1, 2 or 4 corresponds to a width pattern. 
For example $(4,4)$, $(2,4)$ and $(4,2)$ would indicate two adjacent tiles of type 4 and width 4, which cannot be by proposition \ref{vertical_tiles}. 
And since half tiles are actually part of a "full" tile, we cannot have $(1,2,1)$ in any width pattern.

\begin{lemma}\label{tanpoints2}
Let $J=J_{1,n}$ with $n>1$ and associated jigsaw group $\Gamma$. 
Then there is a triple of tangency points  $p_1,p_2,p_3$ of  $\tilde{C}$ such that, 
if $L\neq Id$ is a horizontal translation by less than $\ell(\Gamma)$, then $L(p_1)$, $L(p_2)$ and $L(p_3)$ are no longer tangency points of $\tilde{C}$.
\end{lemma}
\proof 
Let $J = J_{m,n}$ be a jigsaw as in theorem \ref{pseudTHEOREM}, so that every tile in $J$ is vertical. 
Let $E$ be the unique tile of type 4 in $J$ that has two exterior sides. 
We will consider cases depending on the congruency of $n$ modulo 3. 

\begin{figure}[h]
\begin{center}
\includegraphics[scale=0.85]{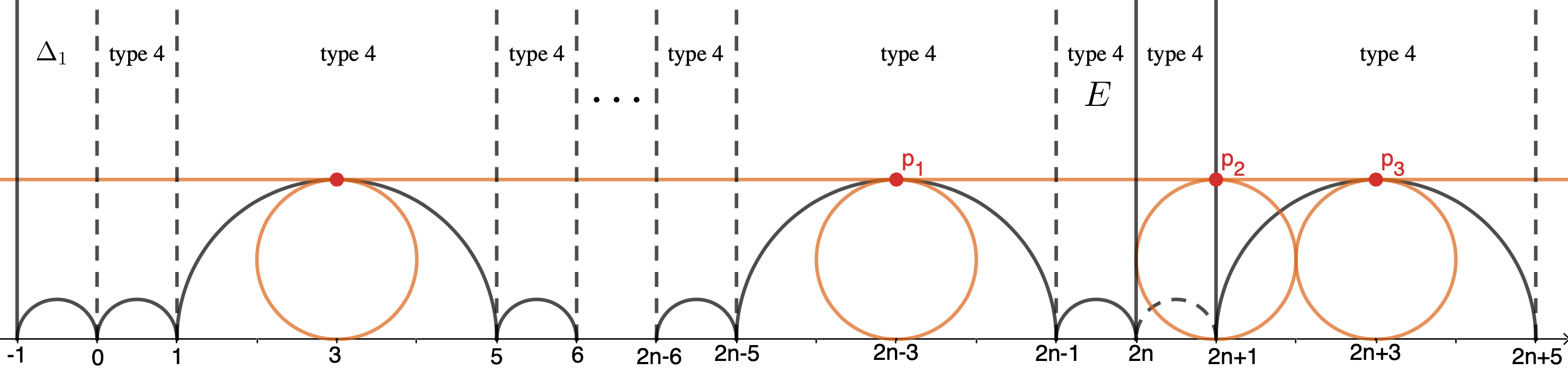} \ \ \ \ \ \ 
\end{center}
\caption{Tangency points for $n\equiv 0 \Mod{3}$}
\label{fig:0mod3tangency}
\end{figure}

$\bullet$ Case 1: $n\equiv 0 \Mod{3}, \ n\geq 3$. 
The tangency points in $\alpha_2\cap J$  are the points $(3+6(j-1),2)$ with $j= 1, \dots, \frac{n}{3}$ (see figure \ref{fig:0mod3tangency}). 
Consider $p_1 = (2n-3,2)$, the last of these points. 
The vertical tile adjacent to $J$ to the right of $[2n,\infty$] is of type 4, and since $n\equiv 0 \Mod{3}$ it has width 1. 
The sides $[\infty,2n]$ and $[2n+1, \infty]$ are exterior of types 1 and 4 respectively, so $p_2 = (2n+1,2)$ is the next tangency point on $\alpha_2$. 
The vertical tile to the right of $[2n+1,\infty$] is of type 4 and width 4. 
Since the side $[2n+1,2n+5]$ is exterior of type 1, the next tangency point on $\alpha_2$ is $p_3 = (2n+3, 2)$. 
Let $L$ be a horizontal translation by less than $\ell(\Gamma)$ and suppose $L(p_i) = p_i' \in \alpha_2$ are tangency points of $\tilde{C}$. 
Since $d_{euc}(p_2',p_3') = d_{euc}(p_2,p_3) = 2$ the width pattern of the tiles between $p_2'$ and $p_3'$ is  $(1,1)$ or $(2)$.

To get the width pattern (1,1) with the desired tangency points $p_2'$ and $p_3'$ we need two tiles $\Delta(\infty, k, k+1)$ and $\Delta(\infty, k+1, k+2)$ of type 4 and width 1 with exterior vertical sides $[k,\infty]$ and $[k+2,\infty]$ of types $\frac{1}{4}$ and 4 respectively. 
Since the tile $E$ of $J$ has interior side of type $\frac{1}{4}$, the tile $\Delta(\infty,k,k+1)$ is not in the orbit $\Gamma\cdot E$. 
Thus the sides $[k,k+1]$ and $[k+1,\infty]$ are interior. 
Then in the adjacent tile $\Delta(\infty, k-4,k)$ the sides $[k-4,\infty]$ and $[k-4,k]$ are interior. 
Since $4 = d_{euc}(p_1, p_2) = d_{euc}(p_1', p_2')$, then $p_1' = (k-4,2)$ which is not a tangency point.

If the width pattern between $p_2'$ and $p_3'$ is (2) we must have a tile $T = \Delta(\infty, k, k+4)$ of type 4 and exterior sides $[\infty,k]$ and $[k,k+4]$. 
The two tiles that follow $T$ to the left must be of type 4 and width 1, with exterior sides  $[k-1,\infty]$ and $[k-2,k-1]$.
Thus the tile $\Delta(\infty, k-2, k-1)$ is a vertical tile in the orbit $\Gamma\cdot E$. 
This implies we obtain $\Delta(\infty, k-2, k-1)$ by translating $E$ by a multiple of $\ell(\Gamma)$ and $L$ must be this translation. 
Therefore $L=Id$.

\begin{figure}[h]
\begin{center}
\includegraphics[scale=0.8]{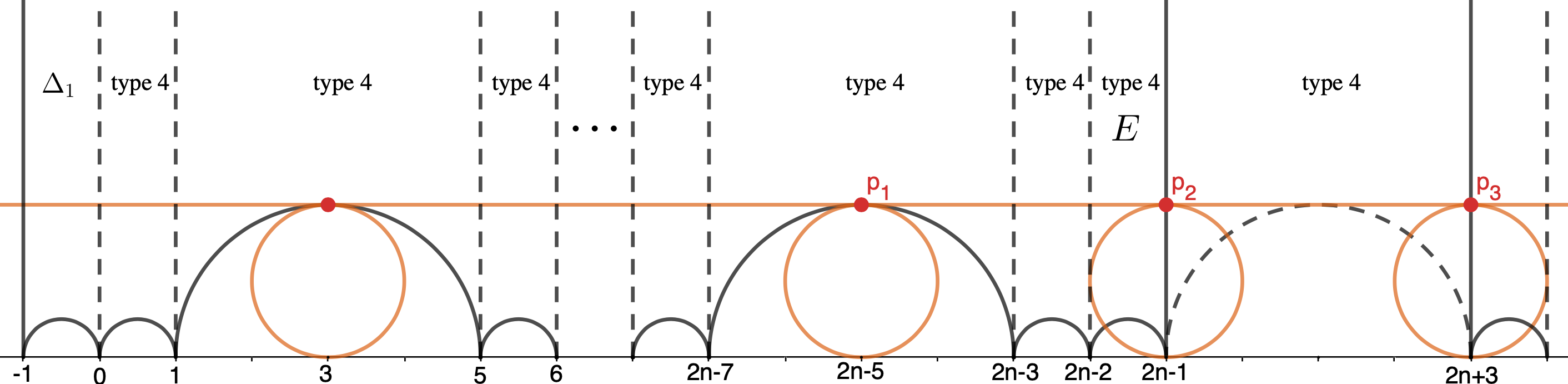} \ \ \ \ \ \ 
\end{center}
\caption{Tangency points for $n\equiv 1\Mod{3}$ for $n\geq 4$}
\label{fig:1mod3tangency}
\end{figure}

$\bullet$ Case 2:  $n \equiv 1\Mod{3}$ and $n\geq 4$. 
In this case the tangency points in $\alpha_2\cap J$ which are not on vertical sides are the points $(3+6(j-1),2)$, with $j= 1\dots \frac{n-1}{3}$. 
Let $p_1=(2n-5,2)$ be the last of such tangency points. 
Since $n\equiv 1 \Mod{3}$ the next vertical tile to the right of $J$ is of type 4 and width 4. 
The sides $[\infty, 2n-1]$ and $[2n+3, \infty]$ are exterior with types 4 and $\frac{1}{4}$ respectively, so $p_2=(2n-1,2)$ and $p_3=(2n+3,2)$ are tangency points of the horocycle (see figure \ref{fig:1mod3tangency}). 
As before, assume $L$ is a horizontal translation by less than $\ell(\Gamma)$ and $p_i'=L(p_i)$ are tangency points. 
Since $d_{euc}(p'_1,p'_2) = d_{euc}(p_1,p_2) =4$ the possible width patterns for tiles between $p_1'$ and $p_2'$ are (4), (1,1,2) , (2,1,1) and (1,1,1,1).

\begin{figure}[h]
\centering
\begin{minipage}{.5\textwidth}
  \centering
  \includegraphics[width=1\linewidth]{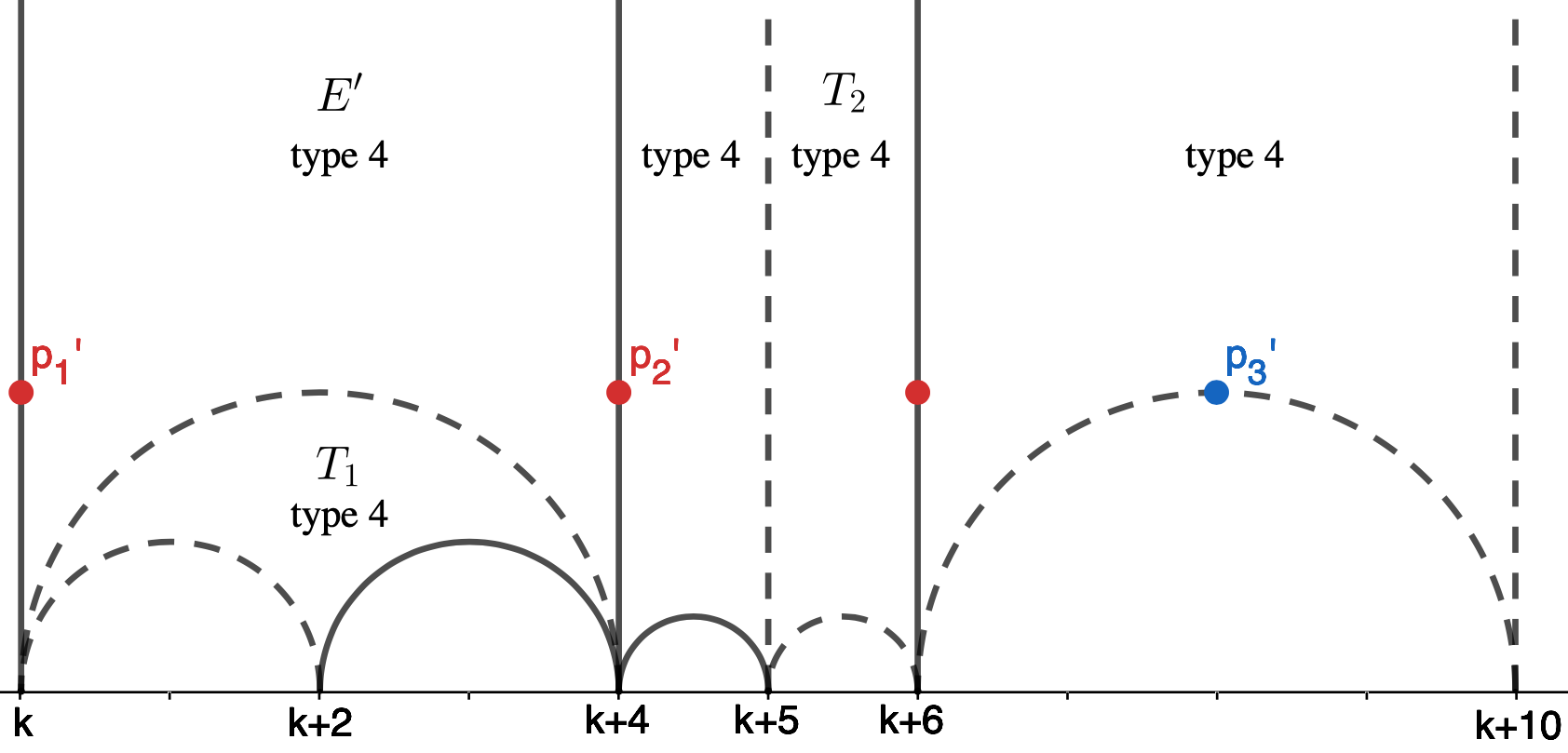}
  \caption{Width pattern (4)}
  \label{fig:case2width4} 
\end{minipage}%
\begin{minipage}{.5\textwidth}
  \centering
  \includegraphics[width=0.7\linewidth]{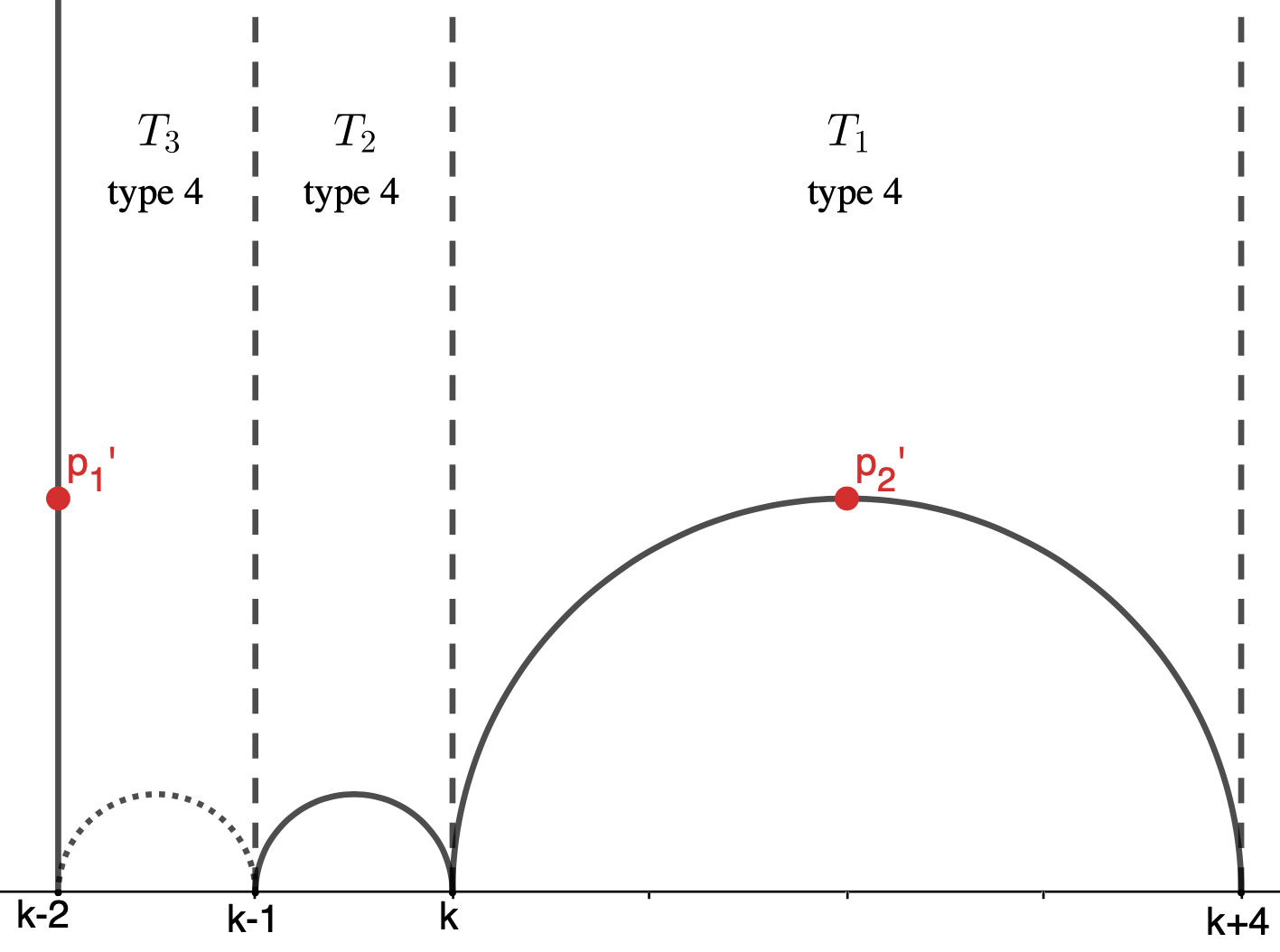}
  \caption{Width pattern (1,1,2)}
  \label{fig:case2width11}
\end{minipage}
\end{figure}

To get the width pattern (4) between $p_1'$ and $p_2'$ we need a tile $E' = \Delta(\infty, k, k+4)$ of type 4 with vertical external sides, so $E'$ is in the $\Gamma$-orbit of $E$ (see figure \ref{fig:case2width4}).
Since $J$ has more than one tile of type 4, the tile $T_1 = \Delta(k,k+4,k+2)$ adjacent to $E'$ must be of type 4 with external side $[k+2,k+4]$. 
Let $R\in \Gamma$ be the $\pi$-rotation about the marked point $p_2'$ on $[k+4,\infty]$. 
Then $T_2 = R(T_1) = \Delta(\infty, k+5, k+6)$ has $[k+6,\infty]$ as an exterior vertical side of type 4. 
It follows that the tile adjacent to $T_2$ to the right has width 4 with interior side $[k+6,k+10]$. 
Since $d_{euc}(p_2', p_3') = d_{euc}(p_2, p_3) = 4 $, then $p_3' = (k+8,2)$ which is not a tangency point.


If the width pattern between $p_1'$ and $p_2'$ is (1,1,2) then we must have a tile $T_1=\Delta(\infty, k, k+4)$ of type 4 where the side $[k,k+4]$ is exterior. Since the tile $E$ in $J$ has interior side of type 1, $T_1 \notin \Gamma\cdot E$. 
Then the sides $[\infty,k]$ and $[k+4,\infty]$ are interior (see figure \ref{fig:case2width11}). 
The tiles $T_2 = \Delta(\infty, k-1, k)$ and $T_3 =\Delta(\infty, k-2, k-1)$ to the left of $T_1$ must be of type 4 and width one, 
with $[k-1,\infty]$ interior of type 1 and $[k-2, \infty]$ exterior of type $\frac{1}{4}$. 
Having two interior vertical sides implies that $T_2$ is in a translation of the initial jigsaw $J$, and thus $T_3$ is too. 
However in $J$ the exterior vertical sides have types 1 and 4, so this tile configuration is not possible. 

If we have width pattern (2,1,1) or (1,1,1,1) between $p_1'$ and $p_2'$ then there is a tile $T_1 = \Delta(\infty, k, k+1)$ of type 4 whose side $[k+1,\infty]$ is exterior of type 4 and has $p_2'$ in it. 
The tile adjacent to $T_1$ to the right has to be of type 4 and width 4, let this be $T_2=\Delta(\infty, k+1, k+5)$. 
Because $d_{euc}(p_2',p_3') = 4$, we have that $p_3' \in [k+5,\infty]$, so this side must be exterior. 
This implies $[k,\infty]$ is exterior too and so $T_3 = \Delta(\infty, k-1,k)$ is of type 4 with exterior sides $[k-1,k]$ and $[k,\infty]$. 
Thus $T_3$ is obtained by translating the tile $E$ in $J$ by a multiple of $\ell(\Gamma)$. 
Since $L$ must be this translation, we get that $L=Id$.

\begin{figure}[h]
\begin{center}
\includegraphics[scale=0.8]{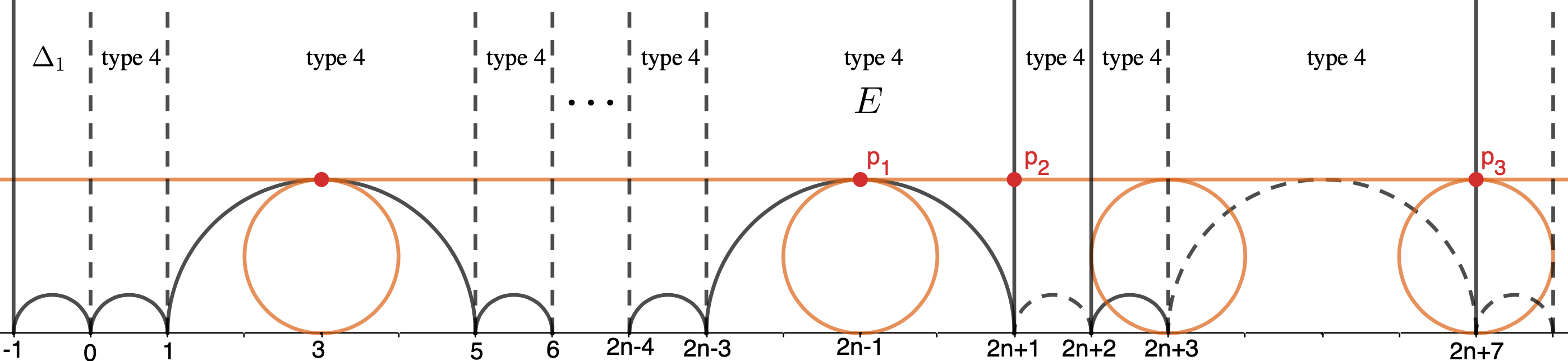} \ \ \ \ \ \ 
\end{center}
\caption{Tangency points for $n\equiv 2 \Mod{3}$}
\label{fig:2mod3}
\end{figure}

$\bullet$ Case 3: $n \equiv 2\Mod{3}$. 
The last tangency points  on $\alpha_2\cap J$ are $p_1 = (2n-1,2)$ and $p_2 = (2n+1,2)$. 
The next tangency point along $\alpha_2$ is $p_3=(2n+7, 2)$ (see figure \ref{fig:2mod3}). 
Let $L$ be a horizontal translation by less than $\ell(\Gamma)$ and suppose $L(p_i) = p_i' \in \alpha_2$ are tangency points of $\tilde{C}$. 
Since $d_{euc}(p_1',p_2') = d_{euc}(p_1,p_2) = 2$ the width pattern of the tiles between $p_1'$ and $p_2'$ is (1,1) or (2).

Suppose the width pattern between $p_1'$ and $p_2'$ is (1,1). 
Then we have two adjacent tiles $\Delta(\infty, k, k+1)$ and  $T_1=\Delta(\infty, k+1, k+2)$ with exterior vertical sides $[k,\infty]$ and $[k+2, \infty]$ of types $\frac{1}{4}$ and 4 respectively. 
Their common side $[k+1,\infty]$ is interior of type 1. 
Since the tile $E$ in $J$ has an exterior sides of type 1, these tiles are not in the $\Gamma$-orbit of $E$. 
Thus the sides $[k,k+1]$ and $[k+1,k+2]$ must be interior.
The tile to the right of $T_1$ has type 4 and width 4, with exterior sides $[k+2,k+6$ and $[k+6,\infty]$. 
Then the tile $T_2 = \Delta(\infty, k+6, k+7)$ must be of type 4 with side $[k+6,\infty]$ of type $\frac{1}{4}$. 
Notice that if $v$ is a vertex of a tile in $J$ which is not $E$ or $\Delta_1$, then three sides of tiles in $J$ meet at $v$, two exterior and one interior. 
Then by looking at the vertex $\infty$ of $T_2$ we get that the side $[k+7, \infty]$ is exterior of type 1. 
Then $\Delta(\infty, k+7,k+8)$ is a tile of type 4 in the tiling and its side $[k+8,\infty]$ is exterior. 
Since $d_{euc}(p_2',p_3') = d_{euc}(p_2,p_3)=6$ we have that $p_3'=(k+8,2)$ is not a tangency point. 

\begin{figure}[h]
\begin{center}
\includegraphics[scale=0.9]{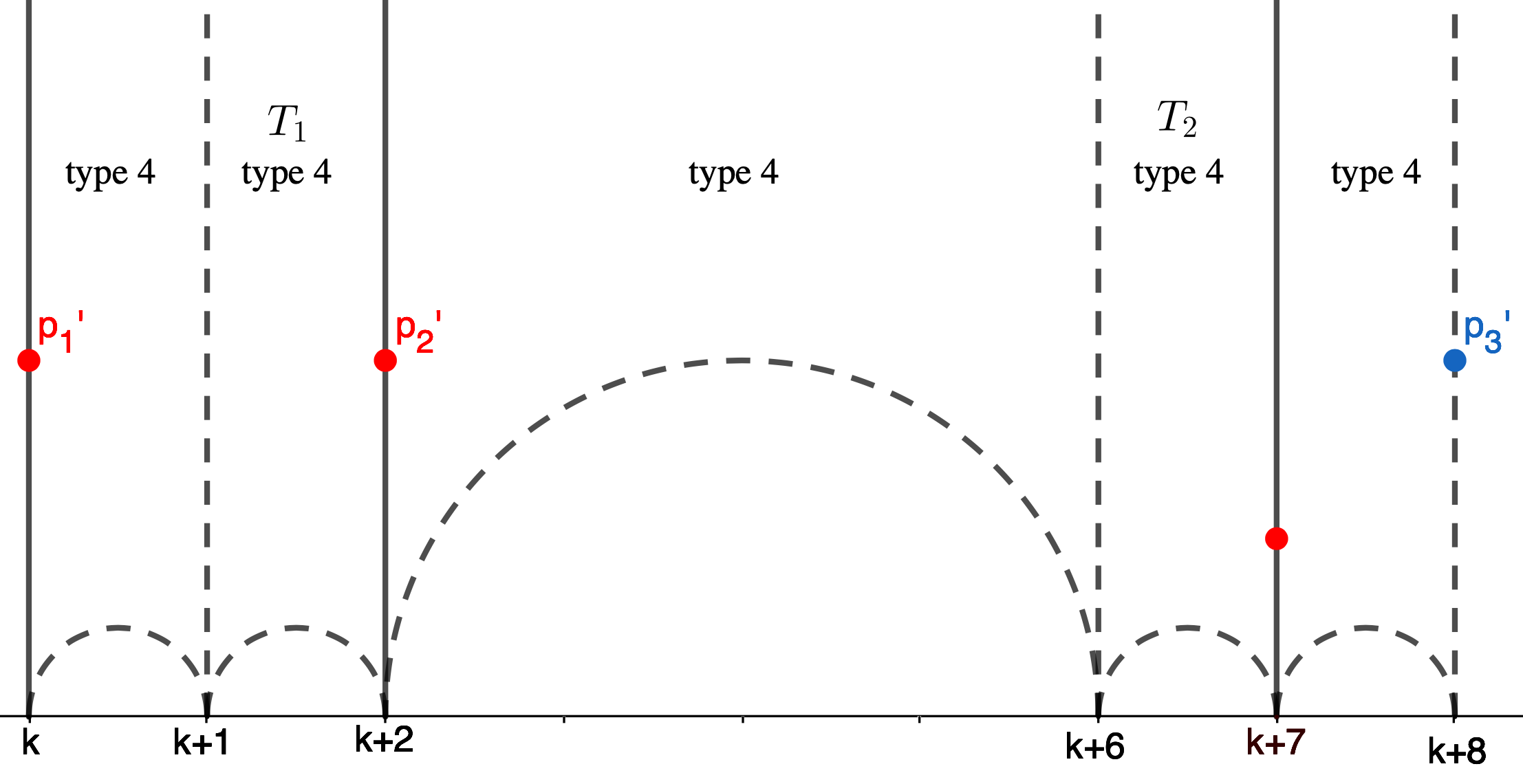} \ \ \ \ \ \ 
\end{center}
\caption{Width pattern (1,1)}
\label{fig:case3width11}
\end{figure}

If (2) is the width pattern between $p'_1$ and $p'_2$ then these tangency points are on a tile $\Delta(\infty, k, k+4)$ with exterior sides $[k, k+4]$ and $[k+4,\infty]$ of types 1 and $\frac{1}{4}$. 
This tile must then be a translation of $E$ in $J$ by a multiple of $\ell(\Gamma)$. 
As before this leads to $L=Id$. 

\vspace{0.5 cm}


\begin{corollary} \label{preservehorocycle}
Let $J_{m,n}$ with $m=1$ or $n=1$ be a jigsaw as in theorem \ref{pseudTHEOREM} and $\Gamma_{m,n}$ its associated jigsaw group. 
Let $\tilde{C}$ be the preimage in $\H^2$ of the maximal horocycle of the orbifold $\mathcal{O} = \H^2/\Gamma_{m,n}$ and $L$ a horizontal translation of $\H^2$ by less than $\ell(\Gamma_{m,n})$. 
Then $L(\tilde{C}) \neq \tilde{C}$. 
\end{corollary}
\proof A horizontal translation that preserves $\tilde{C}$ must preserve the set of tangency points of $\tilde{C}$. Lemmas \ref{tanpoints1} and \ref{tanpoints2} show this is not possible for a translation by less than $\ell(\Gamma_{m,n})$.

\vspace{0.5 cm}


\begin{proposition}\label{non_commens}
Two distinct groups in the families $\Gamma_{1,n}$ and $\Gamma_{m,1}$ are non-commensurable.
\end{proposition}
\proof Let $\Gamma = \Gamma_{m,n}$ be a jigsaw group with $m=1$ or $n=1$ and let $[\Gamma]$ be its commensurability class.
The group $\Gamma$ is non-arithmetic by proposition \ref{nonarith}, so its commensurator $\text{Comm}(\Gamma)$ is the unique maximal element in $[\Gamma]$ \cite{Margulis}. 
In terms of covering spaces $\text{Comm}(\Gamma)$ is the fundamental group of a unique minimal orbifold $\mathcal{O}' =\H^2/\text{Comm}(\Gamma)$ which is finitely covered by any other orbifold $\H^2/G$ with $G\in [\Gamma]$.

Suppose that $\Gamma$ is a proper subgroup of $\text{Comm}(\Gamma)$ and let $\mathcal{O} = \H^2/\Gamma$. 
Let $C$ be the maximal horocycle in $\mathcal{O}'$ and $\tilde{C}$ be the preimage of $C$ to $\H^2$.  
The orbifold $\mathcal{O}$ has a single cusp, so $\tilde{C}$ is the preimage of the maximal horocycle in $\mathcal{O}$ too. 
Since $\mathcal{O}$ covers $\mathcal{O}'$ the lift $\tilde{C}$ must be invariant under a horizontal translation by $k$ where $0<k< \ell(\Gamma)$. 
However, this would contradict corollary \ref{preservehorocycle}. 
Therefore $\Gamma=\text{Comm}(\Gamma)$.

\vspace{0.5 cm}

Propositions \ref{cuspsareQ} and \ref{non_commens} complete the proof of theorem 1. 


\section{Non-pseudomodular Weierstrass groups.}\label{sectionWEIR}
In this final section we prove theorem \ref{WnTHEOREM} which states the Weierstrass groups $W_n$ with $n\geq 6$ congruent to 0, 2 or 6 mod 8 are not pseudomodular. 

By definition a pseudomodular group $\Gamma$ is discrete, so no element in $\R\cup\{\infty\}\equiv \partial_\infty \H^2$ can be simultaneously fixed by a parabolic and a hyperbolic element in $\Gamma$. 
Then to see that a given $\Gamma<PGL(2,\Q)$ is not pseudomodular it suffices to find a hyperbolic element in $\Gamma$ that fixes a rational number. 
Following \cite{LR} we call such a hyperbolic element a \textit{special element} of $\Gamma$, and its fixed points \textit{special points} in $\Q$. 
By constructing special elements we prove infinitely many of the $W_n$ jigsaw groups are not pseudomodular. 
This result provides a partial answer to whether all but finitely many Weierstrass groups are non-pseudomodular, a question posed at the end of \cite{LTV}.

\begin{proposition}
The Weierstrass groups $W_n$ with $n>2$ and congruent to 0, 2 or 6 mod 8 contain a special element.
\end{proposition}
\proof From definition \ref{Weir_groups} the generators of $W_n$ can be represented by the matrices
\begin{equation}\nonumber
	a=\begin{pmatrix}1 & 2 \\ -1 & -1 \end{pmatrix}, \ \ 
	b=\sqrt{n}\begin{pmatrix}1 & 1\\ \frac{-n-1}{n} & -1 \end{pmatrix}, \ \ 
	c= \frac{1}{\sqrt{n}}\begin{pmatrix} 0 & n \\ -1 & 0 \end{pmatrix}
\end{equation}
in $PSL(2,\R)$. 
Recall that a matrix in $PSL(2,\R)$ represents a hyperbolic element of the isometries of $\H^2$ if its trace is bigger than 2. 
Let us examine each congruency class separately. 


$\bullet$ Case 1: $ n = 8k$, for $k\geq 1$.
Consider $A = cba = \begin{pmatrix} 1& 8k+2 \\ 0 & 1 \end{pmatrix} \in W_n$. 
A direct calculation shows that $cabaA^{k-1}caba$ has trace $4k + \frac{1}{4k}>2$ and  its fixed points are the integers $-4k$ and $4k-2$, so this is a special element in $W_{8k}$.


$\bullet$ Case 2: $ n = 8k+2$, for $k\geq 1$.  Let 
	$A = abc = \begin{pmatrix} 1& -8k-4 \\ 0 & 1 \end{pmatrix} \in W_{8k+2}. 	$
It can be directly calculated that  the rationals $\frac{2}{8k+1}$ and $\frac{-8k-2}{4k+3}$ are fixed points of $cA^{4k-1}ababaA^{-k+1}ca$, thus this  is a special element in $W_{8k+2}$.

$\bullet$ Case 3: $n = 8k+6$, for $k\geq 0$. 
Let $ A =  cba = \begin{pmatrix}1 & 8k+8 \\ 0& 1 \end{pmatrix} \in W_{8k+6}$. 
A direct computation shows that  the element $aA^kcababac$ fixes 1 and has trace $9 + 24 k + 16 k^2 + \frac{1}{(3 + 4 k)^2},$ which is greater than 2. 
Thus we have found a special element in $W_{8k+6}$. 

\vspace{0.25 cm}

Theorem 2 immediately follows from the previous proposition. \\

We have found that $W_n$ has a special element for small values of $n\equiv 4 \Mod{8}$, 
though no clear pattern in the fixed points or the word of the  element is clear. 
The computer program we developed to obtain these examples tries to determine whether a given rational is a special point by exploring its $\Gamma$-orbit. 
A survey of whether $W_n$ has a special for $n\leq 28$ follows.



\begin{table}[h]
  \begin{center}
	\begin{tabular}{r |r | r | r }
	\textbf{Group} & \textbf{Pseudomodular or special}&\textbf{Group} & \textbf{Pseudomodular or special}\\
	\hline      
      $W_1$ &  pseudomodular \cite{LTV} 		&$W_{15}$ & contains a special \\
      \hline      
      $W_2$ & pseudomodular \cite{LTV} 		&$W_{16}$ &   contains a special \\
      \hline
      $W_3$ & contains a special 		&$W_{17}$ &  contains a special \\
      \hline
      $W_4$ & pseudomodular (theorem \ref{pseudTHEOREM})			&$W_{18}$ &   contains a special \\
      \hline
      $W_5$ &  contains a special 		&$W_{19}$ &   could not be determined \\
      \hline
      $W_6$ &  contains a special 		&$W_{20}$ &  contains a special \\
      \hline
      $W_7$ &  contains a special 		&$W_{21}$ &  contains a special \\
      \hline
      $W_8$ &  contains a special 		&$W_{22}$ &  contains a special \\
      \hline
      $W_9$ &  contains a special 		&$W_{23}$ &   contains a special \\
      \hline
      $W_{10}$ &  contains a special 		&$W_{24}$ &   contains a special \\
      \hline
      $W_{11}$ & could not be determined 		&$W_{25}$ &   could not be determined \\
      \hline                 
      $W_{12}$ &  contains a special 		&$W_{26}$ &   contains a special \\
      \hline
      $W_{13}$ & could not be determined 		&$W_{27}$ &   contains a special \\
      \hline
      $W_{14}$ &  contains a special 		&$W_{28}$ &   contains a special \\
      \hline                                                       
	\end{tabular}
	
	
	\caption{Survey of small Weierstrass groups}
	\label{tab:table1}

  \end{center}
\end{table}


\end{document}